\documentclass[a4paper,12pt]{article}

\usepackage{amsfonts}
\usepackage{amscd,color}
\usepackage{amsmath,amsfonts,amssymb,amscd}
\usepackage{indentfirst,graphicx,epsfig}
\usepackage{graphicx,psfrag}
\input{epsf}
\usepackage{graphicx}
\usepackage{epstopdf}
\usepackage{caption}

\newtheorem{thm}{Theorem}[section]

\newtheorem{lem}[thm]{Lemma}

\newtheorem{cor}[thm]{Corollary}
\newtheorem{pro}[thm]{Proposition}

\def\qed{\hfill \nopagebreak\rule{5pt}{8pt}}

\title{\textbf{Nordhaus-Gaddum-type theorem for \\
conflict-free connection number of graphs\footnote{Supported by NSFC Nos. 11371205, 11531011, 11601254 and 11551001, and the SFQP Nos. 2016-ZJ-948Q and 2014-ZJ-907), and the project on the
key lab of IOT of Qinghai province (No. 2017-Z-Y21).}}}
\author{\small Hong Chang$^{1}$, \ Zhong Huang$^{1}$, \ Xueliang Li$^{1,2}$, \ Yaping Mao$^{2}$,  \ Haixing Zhao$^{2}$\\[0.2cm]
\small $^{1}$Center for Combinatorics and LPMC \\
\small Nankai University, Tianjin 300071, China\\[0.2cm]
\small $^{2}$School of Mathematics and Statistics \\
\small Qinghai Normal University, Xining, Qinghai 810008, China\\[0.2cm]
{\small Email: changh@mail.nankai.edu.cn, 2120150001@mail.nankai.edu.cn,}\\
{\small lxl@nankai.edu.cn, maoyaping@ymail.com; h.x.zhao@163.com}\\
}
\date{}
\begin{document}
\maketitle
\begin{abstract}
An edge-colored graph $G$ is \emph{conflict-free connected} if,
between each pair of distinct vertices, there exists a path
containing a color used on exactly one of its edges. The
\emph{conflict-free connection number} of a connected graph $G$,
denoted by $cfc(G)$, is defined as the smallest number of colors
that are needed in order to make $G$ conflict-free connected. In
this paper, we determine all trees $T$ of order $n$ for which
$cfc(T)=n-t$, where $t\geq 1$ and $n\geq 2t+2 $. Then we prove that
$1\leq cfc(G)\leq n-1$ for a connected graph $G$, and characterize
the graphs $G$ with $cfc(G)=1,n-4,n-3,n-2,n-1$, respectively.
Finally, we get the Nordhaus-Gaddum-type theorem for the 
conflict-free connection number of graphs, and prove that
if $G$ and $\overline{G}$ are connected, then
$4\leq cfc(G)+cfc(\overline{G})\leq n$ and
$4\leq cfc(G)\cdot cfc(\overline{G})\leq2(n-2)$, and moreover,
$cfc(G)+cfc(\overline{G})=n$ or $cfc(G)\cdot cfc(\overline{G})=2(n-2)$
if and only if one of $G$ and $\overline{G}$ is a tree with maximum
degree $n-2$ or a $P_5$, and the lower bounds are sharp.\\[2mm]
Nordhaus-Gaddum-type.
\end{abstract}

\section{Introduction}

All graphs in this paper are undirected, finite and simple. We
follow \cite{BM} for graph theoretical notation and terminology not
described here. Let $G$ be a graph. We use $V(G), E(G), n(G), m(G)$,
and $\Delta(G)$ to denote the vertex-set, edge-set, number of
vertices, number of edges, and maximum degree of $G$, respectively.
For $v \in V(G)$, let $N(v)$ denote the neighborhood of $v$ in $G$,
and let $d(v)$ denote the degree of $v$ in $G$, and
$d_F(v)$ denote the degree of $v$ in a subgraph $F$ of $G$.
Given two graphs $G$ and $H$, the \emph{union} of $G$ and $H$,
denoted by $G\cup H$, is the graph with vertex-set $V(G)\cup V(H)$
and edge-set $E(G)\cup E(H)$. The \emph{join} of $G$ and $H$,
denoted by $G+H$, is obtained from $G\cup H$ by joining each vertex
of $G$ to every vertex of $H$.

Let $G$ be a nontrivial connected graph with an associated
\emph{edge-coloring} $c : E(G)\rightarrow \{1, 2, \ldots, t\}$, $t
\in \mathbb{N}$, where adjacent edges may have the same color. If
adjacent edges of $G$ are assigned different colors by $c$, then $c$
is a \emph{proper (edge-)coloring}. For a graph $G$, the minimum
number of colors needed in a proper coloring of $G$ is referred to
as the \emph{edge-chromatic number} of $G$ and denoted by
$\chi'(G)$. A path of an edge-colored graph $G$ is said to be a
\emph{rainbow path} if no two edges on the path have the same color.
The graph $G$ is called \emph {rainbow connected} if every pair
of distinct vertices of $G$ is connected by a rainbow path in $G$.
An edge-coloring of a connected graph is a \emph{rainbow connection
coloring} if it makes the graph rainbow connected. This concept of
rainbow connection of graphs was introduced by Chartrand et
al.~\cite{CJMZ} in 2008. For a connected graph $G$, the \emph{rainbow connection number}
$rc(G)$ of $G$ is defined as the smallest number of colors
that are needed in order to make $G$ rainbow connected. The reader
who are interested in this topic can see \cite{LSS, LS} for a survey.

Inspired by rainbow connection coloring and proper coloring in graphs, Andrews
et al.~\cite{ALLZ} and Borozan et al.~\cite{BFGMMMT} introduced the
concept of proper-path coloring. Let $G$ be a nontrivial connected
graph with an edge-coloring. A path in $G$ is called a \emph{proper
path} if no two adjacent edges of the path receive the same color.
An edge-coloring $c$ of a connected graph $G$ is a \emph{proper-path
coloring} if every pair of distinct vertices of $G$ are connected by
a proper path in $G$. And if $k$ colors are used, then $c$ is called
a \emph{proper-path $k$-coloring}. An edge-colored graph $G$ is
\emph{proper connected} if any two vertices of $G$ are connected by
a proper path. For a connected graph $G$, the minimum number of colors
that are needed in order
to make $G$ proper connected is called the \emph{proper connection
number} of $G$, denoted by \emph{$pc(G)$}. Let $G$ be a nontrivial
connected graph of order $n$ and size $m$ (number of edges). Then
we have that $1\leq pc(G) \leq \min\{\chi'(G), rc(G)\}\leq m$.
For more details, we refer to \cite{GLQ,LLZ,LWY} and a dynamic
survey \cite{LM}.

A coloring of vertices of a hypergraph $H$ is a called
\emph{conflicted-free} if each hyperedge $E$ of $H$ has a vertex of
unique color that does not get repeated in $E$. The smallest number
of colors required for such a coloring is called the
\emph{conflict-free chromatic number} of $H$. This parameter was
first introduced by Even et al. \cite{ELR} in a geometric setting,
in connection with frequency assignment problems for cellular
networks. One can find many results on conflict-free coloring, see
\cite{CKP,CT,PT}.

Recently, Czap et al. \cite{CJV} introduced the concept of
conflict-free connection of graphs. An edge-colored graph $G$ is
called \emph{conflict-free connected} if each pair of distinct
vertices is connected by a path which contains at least one color
used on exactly one of its edges. This path is called a
\emph{conflict-free path}, and this coloring is called a
\emph{conflict-free connection coloring} of $G$. The
\emph{conflict-free connection number} (or, {\it cfc number}, for short) of a connected graph $G$,
denoted by $cfc(G)$, is the smallest number of colors needed to
color the edges of $G$ so that $G$ is conflict-free connected. In
\cite{CJV}, they showed that it is easy to get the conflict-free
connection number for $2$-connected graphs and very difficult for
other connected graphs, including trees.

A Nordhaus-Gaddum-type result is a (tight) lower or upper bound on
the sum or product of the values of a parameter for a graph and its
complement. The name ``Nordhaus-Gaddum-type" is given because
Nordhaus and Gaddum \cite{NG} first established the following type
of inequalities for chromatic numbers in 1956. They proved that if
$G$ and $\overline{G}$ are complementary graphs on $n$ vertices whose
chromatic numbers are $\chi(G)$ and $\chi(\overline{G})$, respectively,
then $$2\sqrt{n}\leq \chi(G)+\chi(\overline{G})\leq n+1, \  n \leq \chi(G)\cdot \chi(\overline{G})\leq (\frac{n+1}{2})^2.$$
Since then, the Nordhaus-Gaddum type relations have received wide
attention: diameter \cite{HR}, domination number \cite{HH,SDK}, connectivity \cite{HV}, generalized edge-connectivity \cite{LY}, rainbow connection number
\cite{CLL}, list coloring \cite{DGM}, Wiener index \cite{LWYA} and some other chemical indices \cite{ZW},
and so on. For more results, we refer to a recent survey paper \cite{AH} by Aouchiche and Hansen.

Let us give an overview of the rest of this paper. In Section $2$, we
present some upper bounds for the conflict-free connection number. In
Section $3$,we determine all trees $T$ of order $n$ for which
$cfc(T)=n-t$, where $t\geq1$ and $n\geq2t+2$. In Section $4$, graphs
$G$ with $cfc(G)=1,n-4,n-3,n-2,n-1$ are respectively characterized.
In Section $5$, we get the Nordhaus-Gaddum-type theorem for the conflict-free
connection number of graphs, and prove that if $G$ and $\overline{G}$
are connected, then $4\leq cfc(G)+cfc(\overline{G})\leq n$ and
$4\leq cfc(G)\cdot cfc(\overline{G})\leq2(n-2)$, and moreover,
$cfc(G)+cfc(\overline{G})=n$ or $cfc(G)\cdot cfc(\overline{G})=2(n-2)$
if and only if one of $G$ and $\overline{G}$ is a tree with maximum
degree $n-2$ or a $P_5$, and the lower bounds are sharp.

\section{Preliminaries}

At the very beginning, we state some fundamental results on the
conflict-free connections of graphs, which will be used in the sequel.

\begin{lem}{\upshape \cite{CJV}}\label{lem2-1}
If $P_n$ is a path on $n$ edges, then $cfc(P)=\lceil
\log_2(n+1)\rceil $.
\end{lem}

It is obvious that $cfc(K_{1,n-1})=n-1$ for $n\geq2$. In \cite{CJV} 
the authors obtained the upper and lower bounds of the conflict-free
connection number for trees in terms of the maximum degree $\Delta$.
\begin{lem}{\upshape \cite{CJV}}\label{lem2-2}
If $T$ is a tree on $n$ vertices with maximum degree $\Delta(T)\geq 3$
and diameter $d(T)$, then
$$
\max\{\Delta(T),\log_2(d(T))\}\leq cfc(T)\leq \frac{(\Delta(T)-2)\log_2(n)}{\log_2(\Delta(T))-1}.
$$
\end{lem}

\begin{lem}{\upshape \cite{CJV}}\label{lem2-3}
If $G$ is a noncomlplete $2$-connected graph, then $cfc(G)=2$.
\end{lem}

A \emph{block} of a graph $G$ is a maximal connected subgraph of $G$
that has no cut-vertex. If $G$ is connected and has no cut-vertex,
then $G$ is a block. An edge is a block if and only if it is a
cut-edge, this block is called \emph{trivial}. Then any nontrivial
block is $2$-connected.
\begin{lem}{\upshape \cite{CJV}}\label{lem2-4}
Let $G$ be a connected graph. Then from its every nontrivial block
an edge can be chosen so that the set of all such chosen edges forms
a matching.
\end{lem}

From Lemma \ref{lem2-4}, we can extend the result
of Lemma \ref{lem2-3} to $2$-edge-connected graphs in the following.

\begin{cor}\label{cor2-3}
Let $G$ be a noncomplete $2$-edge-connected graph. Then
$cfc(G)=2$.
\end{cor}
\begin{pf}
Since $G$ is not a complete graph, it follows that $cfc(G)\geq 2$.
It suffices to show that $cfc(G)=2$. From Lemma \ref{lem2-4}, we can
choose an edge in each nontrivial block such that the set $S$ of
such chosen edges forms a matching. Then we color the edges from $S$ with
color $2$ and color the remaining edges of $G$ with color $1$. It is
easy to check that this coloring is a conflict-free connection
coloring of $G$. Thus, $cfc(G)=2$.\qed
\end{pf}

\vskip0.3cm

Let $C(G)$ be the subgraph of $G$ induced on the set of
cut-edges of $G$, and let $h(G)=\max\{cfc(T):$ $T$ is
a component of C(G)$\}$. The following theorem provides
a sufficient condition for graphs $G$ with $cfc(G)=2$.
\begin{lem}\label{lem2-0}\cite{CJV}
If $G$ is a connected graph and $C(G)$ is a linear forest
whose each component has an order $2$, then $cfc(G)=2$.
\end{lem}

\begin{lem}{\upshape \cite{CJV}}\label{lem2-6}
If $G$ is a connected graph, then
$$
h(G)\leq cfc(G)\leq h(G)+1.
$$
Moreover, the bounds are sharp.
\end{lem}

Note that it is supposed to define $h(G)=1$ for the
$2$-edge-connected graph $G$ in addition. Next, we give a sufficient
condition such that the lower bound is sharp in Lemma \ref{lem2-6}
for $h(G)\geq 2$.
\begin{pro}\label{pro2-7}
Let $G$ be a connected graph with $h(G)\geq 2$. If there exists a
unique component $T$ of $C(G)$ such that $cfc(T)=h(G)$, then
$cfc(G)=h(G)$.
\end{pro}
\begin{pf}
Let $T_1,T_2,\ldots,T_s$ be the components of $C(G)$, and let
$B_1,B_2,\ldots,B_r$ be the nontrivial blocks of $G$. Suppose that
$T_1$ is the unique component $T$ of $C(G)$ with $cfc(T)=h(G)$. We
provide an edge-coloring of $G$ as follows. We first color the edges
of $T_1$ with $h(G)$ colors $\{1,\cdots,h(G)\}$ such that $T_1$ is
conflict-free connected. Then color the edges of $T_i$ with at most
$h(G)-1$ colors $\{1,\cdots, h(G)-1\}$ such that $T_i$ is
conflict-free connected for $2\leq i\leq s$. Next, we color the
edges of $B_i$ for $1\leq i \leq r$. By Lemma \ref{lem2-4}, we
choose an edge in $B_i$ such that the set $S$ of such chosen edges
forms a matching. We color the edges from $S$ with color $h(G)$ and
color the remaining edges of $G$ with color $1$. Note that this
coloring is a conflict-free connection coloring of $G$. So, we have
$cfc(G)= h(G)$.\qed
\end{pf}

\vskip0.3cm

Recall that the \emph{edge-connectivity} of a connected graph $G$,
denoted by $\lambda(G)$, is the minimum size of an edge-subset whose
removal from $G$ results a disconnected graph.

The following result will be useful in our discussion.
\begin{lem}\label{lem2-8}
Let $G$ be a connected graph of order $n$ with $\lambda(G)=1$. Then
$cfc(G)\leq n-2r$, where $r$ is the number of nontrivial blocks of
$G$.
\end{lem}
\begin{pf}
Let $B_1,B_2,\ldots,B_r$ be the nontrivial blocks of $G$, and let
$T_1,T_2,\ldots,T_s$ be the components of $C(G)$. Clearly,
$E(G)=E(B_1)\cup E(B_2)\cup \cdots \cup E(B_r)\cup E(T_1)\cup
E(T_2)\cup \cdots \cup E(T_s)$, and $n=|V(G)|=|V(B_1)|+
|V(B_2)|+\cdots+|V(B_r)|+|V(T_1)|+|V(T_2)|+ \cdots
+|V(T_s)|-(r+s-1)$. Let $T_i$ be a component of $C(G)$ with
$h(G)=cfc(T_i)\leq |V(T_i)|-1$. Since $|V(B_i)|\geq3$ for $1\leq i
\leq r$, it follows that
\begin{eqnarray*}
|V(T_i)|&\leq&|V(T_1)|+|V(T_2)|+ \cdots +|V(T_s)|-2(s-1)\\
&=&n-(|V(B_1)|+ |V(B_2)|+ \cdots +|V(B_r)|)\\
&&+(r+s-1)-2(s-1)\leq n-3r+(r+s-1)-2(s-1)\\
&=&n-2r-s+1\leq n-2r,
\end{eqnarray*}
and hence $h(G)=cfc(T_i)\leq |V(T_i)|-1\leq n-2r-1$. Thus we have
$cfc(G)\leq h(G)+1\leq n-2r$.\qed
\end{pf}

\vskip0.3cm

The following result is an immediate consequence of Lemma
\ref{lem2-8}.
\begin{cor}\label{cor2-9}
Let $G$ be a connected graph of order $n$ with $\lambda(G)=1$. If
$G$ has a unique nontrivial block $B$, then $cfc(G)\leq
n+1-|V(B)|$.
\end{cor}

Next, we give some upper bounds for $cfc(G)$, which will be useful
in our discussion.

\vskip0.3cm

It is clear that the addition of an edge to $G$ can
not increase $cfc(G)$.
\begin{pro}\label{pro2-5}
If $G$ is a nontrivial connected graph and $H$ is
a connected spanning subgraph of $G$,
then $cfc(G)\leq cfc(H)$.
\end{pro}

\begin{pro}\label{pro2-10}
If $G$ is a connected graph with size $m_G$ and $H$ is a connected
subgraph of $G$ with size $m_H$, then $cfc(G)\leq cfc(H)+m_G-m_H$.
\end{pro}
\begin{pf} Suppose that $a=cfc(H)$ and $b=cfc(H)+m_G-m_H$.
Let $c_H$ be a conflict-free connection coloring of $H$
using the colors $1,\cdots,a$. Then $c_H$
can be extended to a conflict-free connection coloring $c_G$ of $G$
by assigning the $m_G-m_H$ colors $a,a+1,\cdots, b$ to the
$m_G-m_H$ edges in $E(G)-E(H)$, which implies that $cfc(G)\leq
cfc(H)+m_G-m_H$. \qed
\end{pf}

\section{Trees with given cfc numbers}

In the sequel, let $K_n$, $K_{s,t}$, $P_n$, and $C_n$ denote the
complete bipartite graph of order $s+t$, complete graph of order
$n$, path of order $n$, cycle of order $n$, respectively. Clearly,
a star of order $t+1$ is exactly $K_{1,t}$. In addition, a double star
is a tree with diameter $3$.

The following theorem indicates that when the maximum degree of a
tree is large, we can give the conflict-free connection number
immediately by its maximum degree.
\begin{thm}\label{th3-1}
Let $T$ be a tree of order $n$, and let $t$ be a natural number such that
$t\geq 1$ and $n\geq 2t+2$. Then $cfc(T)=n-t$ if and only if
$\Delta(T)=n-t$.
\end{thm}

We proceed our proof by the following two lemmas.
\begin{lem}\label{lem3-2}
Let $T$ be a tree of order $n$, and let $t$ be a natural number such that
$t\geq 1$ and $n\geq 2t+2$. If $cfc(T)=n-t$, then $\Delta(T)=n-t$.
\end{lem}
\begin{pf}
Let $T$ be a tree of order $n$ such that $cfc(T)=n-t$ for $1\leq
t\leq 2n+2$. It suffices to show that $\Delta(T)=n-t$. Let $T'$ be a
subtree of $T$ obtained from $T$ by deleting all pendant vertices.
Set $E(T')=\{e_1,\ldots, e_r\}$. Consider any edge $e_i$ of $E(T')$.
There exist two components $A_i$ and $B_i$ of $T-e_i$ for $1\leq i
\leq r$. Without loss of generality, assume that $e(A_i)\geq
e(B_i)\geq 1$ for each $i \ (1\leq i \leq r)$. Let $B_m$ ($1\leq m
\leq r$) be the maximum component in $\{B_1,\ldots, B_r\}$.

\textbf{Claim 1.} $e(B_m)\leq t-1$.

\noindent\textbf{Proof of Claim 1:} Assume, to the contrary, that
$e(B_m)\geq t$.  We define an edge-coloring of $T$ as follows: color
the edges of $A_i$ and $e_i$ with distinct colors, then color the
edges of $B_i$ with distinct colors that assigned to the edges of
$A_i$. Note that $A_i$ and $B_i$ are conflict-free connected, and
the color assigned to $e_i$ is used only once. Thus, it is easy to
see that this coloring is a conflict-free connection coloring, and
hence $cfc(T)\leq n-1-t$, a contradiction. \qed.

\textbf{Claim 2.} For each pair of $B_p$ and $B_q$,  $V(B_p)\cap V(
B_q)=\emptyset$ or $V(B_p)\subseteq V(B_q)$ or $V(B_q)\subseteq V(
B_p)$.

\noindent\textbf{Proof of Claim 2:} Note that $T-e_p-e_q$ has three
components, say $X$, $Y$, and $Z$. If there exists one of $B_p,B_q$
such that it contains two components of $T-e_p-e_q$, then without
loss of generality, we assume $B_q=X\cup e_p \cup Y$. Then $A_q=Z$
and $e(Z)\geq e(B_q)$, and hence $\{A_p,B_p\}=\{X,Y\cup e_q \cup
Z\}$. Since $e (A_p)\geq e(B_p) $ and $e(Y\cup e_q \cup Z)> e(Z)\geq
e(B_q)>e(X)$, it follows that $X=B_p$, and hence $V(B_p)\subseteq
V(B_q)$. If both $B_p$ and $B_q$ have the property that each of them
contains only one component of $T-e_p-e_q$, then $B_p\cap
B_q=\emptyset$.\qed

Let $H=\bigcup _{1\leq i\leq r}B_i$ be a subgraph of $T$. It follows
from Claim $2$ that $B_m$ is also a maximum component of $H$.

Let $F$ be a subgraph obtained from $T$ by deleting the edges of
$H$, and then deleting isolated vertices. Then we have the following
claim.

\textbf{Claim 3.}  $F$ is a star.

\noindent\textbf{Proof of Claim 3:} Assume, to the contrary, that
$F$ is not a star. From Claim $2$, $F$ is connected. Since $F$ is
not a star, it follows that $F$ contains some edge $e_k$ such that
$e_k\in E(T')$. Since $T-e_k$ has two components $A_k$ and $B_k$, it
follows that some edges of $F$ are contained in $B_k$, which is
impossible. Thus, $F$ is a star.\qed

\textbf{Claim 4.} $e(H)=t-1$.

\noindent\textbf{Proof of Claim 4:} At first, we show that $e(H)\leq
t-1$. By contradiction, assume that $e(H)\geq t$. Since $e(H)\geq
t$, one may take $t$ edges of $H$ such that the edges of $B_m$ must
be chosen. Let $D$ be the subgraph of $H$ induced on the set of
these $t$ edges, and let $D_1,\ldots, D_d$ be the components of $D$, in
which $D_1=B_m$. Clearly, $D_1$ is a maximum component of $D$. Note
that for each $D_i$ ($2\leq i \leq d$), there exists one $B_{a_i}$
such that $D_i \subseteq B_{a_i}$. Since $t=e(D)=e(D_1)+\cdots
+e(D_d)$ and $e(D_i)\geq1$ for $1\leq i\leq d$, it follows that
$e(B_m)=e(D_1)\leq t- (d-1)=t-d+1$, and hence $e(B_m)=e(D_1)\leq
n-1-t-d$ for $n\geq 2t+2$. Next, we will provide a conflict-free
connection coloring of $T$ with at most $n-t-1$ colors.

\textbf{Step 1.} Color the edges of
$T\setminus(D\cup\{e_{a_1}\cup\ldots e_{a_d}\})$ with $n-1-t-d$
distinct colors $1, \cdots, n-1-t-d$, and color the edges of
$\{e_{a_1},\ldots, e_{a_d}\}$ with $d$ fresh colors $n-t-d,
\cdots, n-1-t$.

\textbf{Step 2.} Color the edges of $D=D_1\cup\ldots\cup D_d$ with
used colors. Color each edge of $D_1$ with distinct colors that are
assigned to the edges in $T\setminus(D\cup\{e_{a_1}\cup\ldots
e_{a_d}\})$. For any other component $D_i$ ($2\leq i \leq d$) of
$D$, $e(D_i)=e(B_{a_i})-e(B_{a_i}\setminus D_i)\leq
e(B_m)-e(B_{a_i}\setminus D_i)\leq n-1-t-d-e(B_{a_i}\setminus D_i)$.
Next, color each edge of $D_i$ with distinct colors that are
assigned to the edges in $T\setminus(D\cup\{e_{a_1}\cup\ldots
e_{a_d}\}\cup (B_{a_i}\setminus D_i))$.

 In order to prove that this coloring is a conflict-free connection coloring,
 it suffices to show that for each pair of vertices $x,y$, there exists a
 conflict-free path between them. Note that $D_i$ is conflict-free connected,
 the edges of $E(T)\setminus E(D)$ are colored
 with distinct colors, and the color assigned to each edge of $\{e_{a_1}\cup\ldots e_{a_d}\}$
 is used only once under this coloring. Thus, we only need to consider the case
 that $x,y$ are in distinct components of $D$. Without loss of generality, assume
 that $x\in D_i$ and $y\in D_j$, where $1\leq i\neq j\leq d$. The edge assigned
 the unique color on the conflict-free path between them is $e_{a_i}$ (or $e_{a_j}$ ).
 Thus, $cfc(T)=n-t-1$, a contradiction. Thus, $e(H)\leq t-1$. By Claim 2,
 we have $F=K_{1,n-1-e(H)}$. If $e(H)\leq t-2$, then
 $\Delta(T)\geq \Delta(F)=n-1-e(H)\geq n-t+1$, and hence $cfc(T)\geq \Delta(T) \geq n-t+1$
 by Lemma \ref{lem2-2}, which is impossible. Therefore, $e(H)=t-1$. \qed

Let $v$ be a non-leaf vertex of $F$ with $d_F(v)=n-1-e(H)=n-t$.
Since $n\geq 2t+2$, it follows that $d_F(v)=n-t \geq  t+2 \geq 3$.
We claim that $d(v)=n-t$. Assume, to the contrary, that that there
exists an edge $f$ incident with $v$ such that $f\in E(H)$. Then
there exists some $B'$ satisfying $f\in E(B')$, and hence the edges
of $F$ are contained in $B'$, which is impossible. Thus,
$d(v)=n-t\geq t+2 $. Consider any leaf $u$ of $F$. Since $u$ is not
adjacent to the other leaves of $F$, we have $d(u)\leq
n-1-(n-t-1)=t< d(v)$. Consider any vertex $w$ of $V(T)\setminus
V(F)$. It is adjacent to at most one vertex of $F$; otherwise, there
exists a cycle in $T$, a contradiction. Thus, $d(u)\leq
n-1-(n-t)=t-1< d(v)$. As a result, $\Delta(T)=d(v)=n-t$.\qed
\end{pf}

\begin{lem}\label{lem3-3}
Let $T$ be a tree of order $n$, and let $t$ be a natural number such that
$t\geq 1$ and $ n\geq 2t+2$. If $\Delta(T)=n-t$, then $cfc(T)=n-t$.
\end{lem}
\begin{pf}
For convenience, we still use the notation in Lemma \ref{lem3-2}.
Since $\Delta(T)=n-t$, then it follows from Lemma \ref{lem2-2} that
$cfc(T)\geq n-t$. Let $cfc(T)= n-t+u=n-(t-u)$ where $0\leq u\leq
t-1$. It is sufficient to show that $u=0$. By Lemma \ref{lem3-2}, we
have that $e(H)=t-u-1$ and $F=K_{1,n-t+u}$. Note that $\Delta(F)=n-t+u$
and so $\Delta(T)\geq \Delta(F)$, but $\Delta(T)=n-t$. Thus, $u=0$.
The proof is complete.\qed
\end{pf}

\section{Graphs with given small or large cfc numbers}

We first give sharp lower and upper bounds of $cfc(G)$ for a
connected graph $G$.
\begin{pro}\label{pro2-11}
Let $G$ be a connected graph of order $n$. Then
$$
1\leq cfc(G)\leq n-1.
$$
\end{pro}
\begin{pf}
The lower bound is trivial. For the upper bound, we assign distinct
colors to the edges of a given spanning tree of $G$, and color the
remaining edges with one used colors. Since this coloring is a
conflict-free connection coloring of $G$, it follows that
$cfc(G)\leq n-1$.\qed
\end{pf}

\vskip 0.3cm

Graphs with $cfc(G)=1$ can be easily characterized.
\begin{pro}\label{pro3-1}
If $G$ is a connected graph of order $n$, then $cfc(G)=1$
if and only if $G=K_n$.
\end{pro}
\begin{pf}
It is obvious that $cfc(K_n)=1$. Conversely, we let $G$ be a
connected graph with $cfc(G)=1$. If $diam(G)\geq 2$, then we let
$x,y$ be two vertices with $d(x,y)=diam(G)$. Since the conflict-free
path between $x$ and $y$ needs at least two colors, it
follows that $cfc(G)\geq 2$, a contradiction. Thus, $diam(G)=1$,
which implies $G=K_n$.\qed
\end{pf}

\vskip 0.3cm

Next, we present a sufficient and necessary condition for a graph $G$ with $cfc(G)=2$ under the case $diam(\overline{G})\geq3$.

\begin{lem}\label{lem5-6}\cite{ZW}
Let $G$ be a connected graph with connected complement $\overline{G}$. Then
\end{lem}
\begin{itemize}
  \item[] $(i)$ if $diam(G)>3$, then $diam(\overline{G})=2$,
  \item[] $(ii)$ if $diam(G)=3$, then $\overline{G}$ has a spanning subgraph which is a double star.
\end{itemize}

\begin{thm}\label{th5-7}
Let $G$ be a connected noncomplete graph such that $diam(\overline{G})\geq3$.
Then $cfc(G)=2$ if and only if there exist at most two cut-edges incident
with any vertex of $G$.
\end{thm}

\begin{pf}
Suppose $cfc(G)=2$. Assume, to the contrary, that there
exist at least three cut-edges incident with some vertex of $G$.
In order to make $G$ conflict-free connected, these cut-edges
need to be assigned three distinct colors. Thus, $cfc(G)\geq 3$,
a contradiction. For the converse, if $G$ is $2$-edge-connected, then $cfc(G)=2$ by Corollary \ref{cor2-3}. Next, we only consider the case that $G$ has at least one cut-edge. If $n=4,5$, then it is easy to see that the result holds. Next, assume that $n\geq6$. In order to complete our proof, we distinguish the following two cases.

\textbf{Case 1.} $diam(\overline{G})>3$.

It follows from Lemma \ref{lem5-6} that $diam(G)=2$. Take a vertex
$v$ of maximum degree in $G$, let $N_1(v)=\{u_1,\ldots,u_{a}\}$
denote the neighborhood of $v$, and
$N_2(v)=V\setminus \{N_1(v)\cup{v}\}=\{w_1,\ldots,w_{b}\}$, 
where $a\geq1$ and $a+b=n-1$. Consider any vertex
$w_i$ ($1\leq i\leq b $) in $N_2(v)$. Suppose $d(w_i)=1$, and 
let $u_{p_i}$ be the unique neighbour vertex. Since $diam(G)=2$,
it follows that $u_{p_i}$ is adjacent to all other vertices
in $N_1(v)$. Let $W=\{w_i| d(w_i)=1\}$, and $U=\{u_{p_i}|u_{p_i}$ is the unique neighbor of $w_i\}$. Since $diam(G)=2$, it follows that $|W|\leq2$ and $|U|\leq 1$. Note that each vertex of $N_2(v)\setminus W$ is adjacent to at least two vertices of $N_1(v)$. If $W=\emptyset$, then it is easy to see that $cfc(G)=2$. Suppose $W\neq\emptyset$. It follows that $G[V\setminus W]$ is $2$-edge-connected. Let $c'$ be the conflict-free connection coloring of $G[V\setminus W]$ with two colors 1 and 2. Then $c'$ can be extended to a conflict-free connection coloring $c$ of $G$ by assigning distinct colors to cut-edges incident with the vertex in $U$ (if these edges exist). It is easy to see that $G$ is conflict-free connected under this coloring $c$, and thus $cfc(G)=2$.

\textbf{Case 2.} $diam(\overline{G})=3$.

It follows from Lemma \ref{lem5-6} that $G$ has a spanning subgraph $T$ which is a double star. We first present the following claim.

\textbf{Claim } If $e$ is a cut-edge of $G$, then $e\in E(T)$.

\noindent\textbf{Proof of Claim:} Assume, to the contrary, that
$e$ is a cut-edge of $G$ and $e\notin E(T)$. Since $T$ is a
spanning tree of $G-e$, it follows that $G-e$ is connected,
which is impossible.\qed

Let $u,v$ be two non-leaf vertices of $T$, and let $A=\{u\}\cup (N_T(u)\setminus \{v\})=\{u,u_1,\ldots, u_a\}$ and $B= \{v\}\cup (N_T(v)\setminus \{u\})=\{v,v_1,\ldots, v_b\}$, where $a+b=n-2$. Let $H_1=G[A]$ and $H_2=G[B]$.

\textbf{Case 2.1.} $uv$ is a cut-edge of $G$.

It follows that each of $u,v$ is incident with at most one cut-edge in $H_1$,
$H_2$, respectively. Note that $C(G)$ has a unique path whose
order is at most $4$, it follows from Lemma \ref{lem2-0} and Proposition \ref{pro2-7} that $cfc(G)=2$.

\textbf{Case 2.2.} $uv$ is not a cut-edge of $G$.

It follows that there exists an edge $e$ other than $uv$ in $E[A,B]$. Since each of $u,v$ is incident with at most two cut-edges in $H_1$, $H_2$, respectively. Let $H_i^{'}$ be the resulting graph obtained by deleting pendent vertices from $H_i$ (if these vertices exists) for $i=1,2$. Note that $H_i^{'}$ is $2$-edge-connected. Then $cfc(H_i)=h(H_i)=2$ by Proposition \ref{pro2-7}. Let $c_i$ be the conflict-free connection coloring of $H_i$ with two colors 1 and 2 for $i=1,2$. Then $c_1$ and $c_2$ can be extended to a conflict-free connection coloring $c$ of $G$ by coloring the edge $uv$ with color $1$, the edge $e$ with color $2$. It is easy to see that $G$ is conflict-free connected under this coloring $c$, and thus $cfc(G)=2$.\qed
\end{pf}

\begin{thm}\label{th4-1}
Let $G$ be a connected graph of order $n \ (n\geq 2)$. Then $cfc(G)=n-1$
if and only if $G=K_{1,n-1}$.
\end{thm}
\begin{pf}
The sufficiency is trivial, and so we only give the proof of the
necessity. Suppose that $cfc(G)=n-1$. We claim that $G$ is a tree.
Assume, to the contrary, that $G$ is not a tree. Then $G$ contains a
cycle, and so there exists a nontrivial block containing this cycle.
By Lemma \ref{lem2-8}, we have $cfc(G)\leq n-2$, which is
impossible. In order to complete our proof, it is sufficient to show
that $diam(G)=2$. If this is not the case, then $diam(G)\geq 3$. Let
$P$ be the path of length $diam(G)$. From Lemma \ref{lem2-1},
$cfc(P)=\lceil \log_2(diam(G)+1)\rceil\leq diam(G)-1$. From
Proposition \ref{pro2-10}, we have $cfc(G)\leq cfc(P)+n-1-diam(G)
\leq n-2$, a contradiction.\qed
\end{pf}

\vskip 0.3cm

For a nontrivial graph $G$ for which $G+uv\cong G+xy$ for every two
pairs $\{u,v\}$ and $\{x,y\}$ of nonadjacent vertices of $G$, the
graph $G+e$ is obtained from $G$ by adding the edge $e$ joining two
nonadjacent vertices of $G$.
\begin{thm}\label{th4-2}
Let $G$ be a connected graph of order $n \ (n\geq 3)$. Then
$cfc(G)=n-2$ if and only if $G$ is a tree with $\Delta(G)=n-2$ for
$n\geq4$, or $G\in\{K_3, K_{1,3}+e, K_{2,2}, K_{2,2}+e, P_5\}$.
\end{thm}
\begin{pf}
If $G$ is a tree with $\Delta(G)=n-2$, then it follows from Lemma
\ref{lem2-2} and Theorem \ref{th4-1} that $cfc(G)=n-2$. From Lemmas
\ref{lem2-1} and \ref{lem2-3} and Proposition \ref{pro3-1}, we have
$cfc(G)=n-2$ if $G\in\{K_3, K_{1,3}+e, K_{2,2},K_{2,2}+e, P_5\}$.
Thus, it remains to verify the converse. Let $G$ be a connected
graph with $cfc(G)=n-2$. If $3 \leq n\leq 5$, then it is easy to
verify that $G$ is a tree of order $5$ with $\Delta(G)=3$, or
$G\in\{ K_3, P_4, K_{1,3}+e, K_{2,2},K_{2,2}+e, P_5\}$. From now on,
we assume $n\geq 6$. In order to prove our result, we present the
following claim.

\textbf{Claim.} $G$ is a tree.

\noindent\textbf{Proof of Claim:} Assume, to the contrary, that $G$
is not a tree. Then $G$ contains a cycle, and so there exists a
nontrivial block containing this cycle. If $G$ has at least two
nontrivial blocks, then it follows from Lemma \ref{lem2-8} that
$cfc(G)\leq n-4$, which is impossible. Suppose that there exists
only one nontrivial block $B$. Let $T_1,T_2,\ldots,T_s$ be the
components of $C(G)$. Note that $E(G)=E(B)\cup E(T_1)\cup E(T_2)\cup
\cdots \cup E(T_s)$ and $n=|V(G)|=|V(B)|+|V(T_1)|+|V(T_2)|+ \cdots
+|V(T_s)|-(r+s-1)$. If $|V(B)|\geq 4$, then it follows from
Corollary \ref{cor2-9} that $cfc(G)\leq n+1-|V(B)|\leq n-3$, a
contradiction. Suppose $|V(B)|=3$. We assign $|V(T_1)|-1+|V(T_2)|-1+
\cdots +|V(T_s)|-1$ distinct colors $\{1, \cdots,
|V(T_1)|-1+|V(T_2)|-1+ \cdots +|V(T_s)|-1 \}$ to the edges of
$T_1\cup T_2\cup\cdots\cup T_s$, and then color the edges of $B$
with three used colors $\{1,2,3\}$. It is easy to check out that
this coloring is a conflict-free connection coloring of $G$, and so
$cfc(G)\leq |V(B)|-3+|V(T_1)|-1+|V(T_2)|-1+ \cdots +|V(T_s)|-1=
n-3$, which is impossible.\qed

Since $n\geq 6$, it follows from Lemma \ref{lem3-2} that
$\Delta(G)=n-2$. We complete the proof.\qed
\end{pf}

\vskip0.3cm

A graph is \emph{unicyclic} if it is connected and contains exactly
one cycle. Note that $K_{1,n-1}+e$ is unicyclic, whose cycle is a
triangle. Next, we consider another class of unicyclic graphs, whose
cycles are also a triangle. Let $S_{a,n-a}$ be a tree with diameter
$3$, such that the two non-leaf vertices have degree $a$ and $n-a$, the
unicyclic graph $U_n$ is obtained from $S_{3,n-3}$ by adding an edge
joining the two neighbouring leaves of the vertex of degree $3$.
Next, we study the conflict-free connection numbers of $K_{1,n-1}+e$ or
$U_n$.

\begin{lem}\label{lem4-3}
If $G$ is $K_{1,n-1}+e$ or $U_n$ with $n\geq5$, then $cfc(G)=n-3$.
\end{lem}
\begin{pf}
If $G$ is $K_{1,n-1}+e$ or $U_n$ with $n\geq5$, then $G$ has the
property that $h(G)\geq 2$ and there exists a unique component $T$ of
$C(G)$. Noticing that $T=K_{1,n-3}$, then $cfc(G)=h(G)=n-3$ by
Proposition \ref{pro2-7} and Theorem \ref{th4-1}.\qed
\end{pf}

\begin{thm}\label{th4-6}
Let $G$ be a connected graph of order $n\geq 4$. Then $cfc(G)=n-3$
if and only if $G$ satisfies one of the following conditions.
\begin{itemize}
\item[]$(i)$ $G$ is a tree with $\Delta(G)=n-3$, where $ n\geq 6$,

\item[]$(ii)$ $G=K_{1,n-1}+e$, where $ n\geq5$,

\item[]$(iii)$ $G=U_n$, where $ n\geq 5$,

\item[]$(iv)$ $G$ is a $2$-edge-connected and non-complete graph of order
$5$,

\item[] $(v)$ $G\in \{K_4, P_6, G_1, G_2, G_3, G_4, G_5,G_6\}$,
where $G_1, G_2, G_3, G_4, G_5,G_6$ are showed in Fig. 1.
\end{itemize}
\begin{figure}[h!]
\centering
\includegraphics[width=1.0\textwidth]{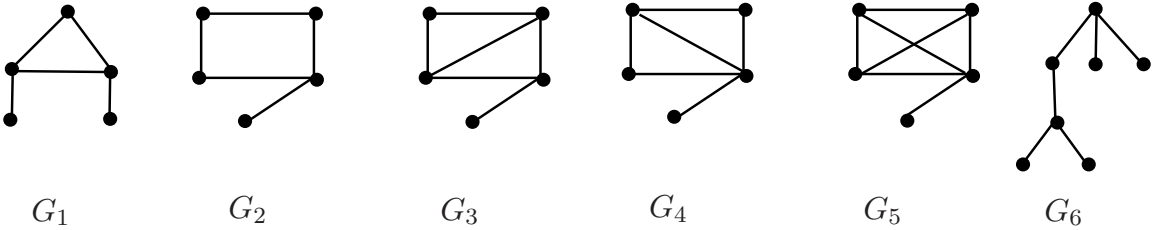}
\caption{Six graphs in Theorem \ref{th4-6} }
\end{figure}
\end{thm}
\begin{pf}
If $G$ is a tree with $\Delta(G)=n-3$ for $n\geq 6$, then $cfc(G)=
n-3$ by Lemma \ref{lem2-2}, Theorems \ref{th4-1} and \ref{th4-2}. If
$G=K_{1,n-1}+e$ or $U_n$ with $n\geq5$, then $cfc(G)=n-3$ by Lemma
\ref{lem4-3}. If $G$ is a $2$-edge-connected and noncomplete graph
of order $5$, then it follows from Corollary \ref{cor2-3} that
$cfc(G)=n-3=2$. Clearly, $cfc(K_4)=cfc(P_6)=n-3$, and each graph
$G_i$ in Fig. 1 satisfies $cfc(G_i)=n-3$ for $1\leq i \leq 6$.
Conversely, let $G$ be a connected graph of order $n\geq 4$ such
that $cfc(G)=n-3$. If $n=4,5$, then $G$ is a $2$-edge-connected and
noncomplete graph of order $5$ or $G\in \{K_4, K_{1,4}+e, U_5,
G_1,G_2, G_3, G_4,G_5 \}$. From now on, we assume $n\geq 6$. We
distinguish the following two cases to show this theorem.

\begin{figure}[h!]
\centering
\includegraphics[width=0.5\textwidth]{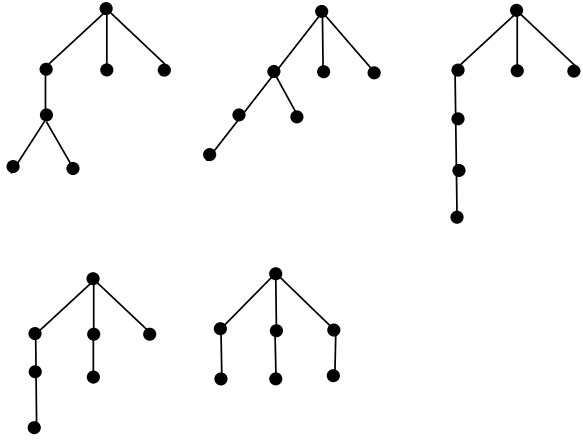}
\caption{Five trees of order $7$ with maximum degree $3$ }
\end{figure}

\textbf{Case 1.} $G$ is a tree of order $n$ with $cfc(G)=n-3$.

For $n=6$, if $\Delta(G)=2$, then $G=P_6$ by Lemma \ref{lem2-1}; if
$\Delta(G)=3$, the result holds trivially; if $\Delta(G)\geq 4$,
then it follows from Lemma \ref{lem2-2} that $cfc(G)\geq 4>n-3$,
which is impossible. Suppose $n=7$. If $\Delta(G)=2$, then $G=P_7$,
and so $cfc(G)=3$ by Lemma \ref{lem2-1}, a contradiction. If
$\Delta(G)=3$, then $G$ is one of five trees in Fig. 2. Note that
the first tree in Fig. 2 has conflict-free connection number $4$,
and the others have conflict-free connection number $3$. Thus, $G$
is the first tree in Fig. 2, which is exactly $G_6$ in Fig. 1. If
$\Delta(G)=4$, then the result always follows. If $\Delta(G)\geq 5$,
then $cfc(G)\geq 5>n-3$ by Lemma \ref{lem2-2}, which is again
impossible. If $n\geq 8$, then $\Delta(G)=n-3$ by Lemma
\ref{lem3-2}.

\textbf{Case 2.} $G$ contains a cycle.

Since $G$ contains a cycle, it follows that there exists a
nontrivial block containing this cycle. Let $T_1,T_2,\ldots,T_s$ be
the components of $C(G)$ such that $h(G)=cfc(T_1)\geq cfc(T_2)\geq
\cdots \geq cfc(T_s)$. If $G$ has at least two
nontrivial blocks, then it follows from Lemma \ref{lem2-8} that
$cfc(G)\leq n-4$, which is impossible. Thus, $G$ has only one
nontrivial block $B$. If $|V(B)|\geq 5$, then $cfc(G)\leq n+1-|V(B)|\leq n-4$
by Corollary \ref{cor2-9}, a contradiction. Suppose $|V(B)|=4$. It
follows from the proof of Lemma \ref{lem2-8} that $cfc(G)
\leq h(G)+1\leq |V(T_1)| \leq
n-|V(B)|+(1+s-1)-2(s-1)\leq n-|V(B)|-s+2 \leq n-2-s$.
If $s\geq 2$, then $cfc(G)\leq n-4$, which is impossible.
If $s=1$, then $cfc(T_1)\leq |V(T_1)|-1 \leq
n-3-1=n-4$, and hence $cfc(G)=h(G)=cfc(T_1)\leq n-4$ by Proposition
\ref{pro2-7}, which is again impossible. Suppose $|V(B)|=3$. It
follows from the proof of Lemma \ref{lem2-8} that $cfc(G)
\leq h(G)+1\leq |V(T_1)| \leq
n-|V(B)|+(1+s-1)-2(s-1)\leq n-|V(B)|-s+2 \leq n-1-s$.
If $s\geq 3$, then $cfc(G)\leq n-4$,
which is impossible. Suppose $s=2$. If $|V(T_2)|\geq 3$, then
$h(G)=cfc(T_1)\leq |V(T_1)|-1 \leq
n-|V(B)|-|V(T_2)|+2-1=n-5$, and so $cfc(G)\leq
h(G)+1\leq n-4$ by Lemma \ref{lem2-6}, a contradiction. If
$|V(T_2)|=2$, then $|V(T_1) |\geq 3$ and so there exists only one
component $T_1$ of $C(G)$ such that $h(G)=cfc(T_1)\leq|V(T_1)|-1\leq
n-4$. Thus, $cfc(G)=h(G)\leq n-4$ by Proposition \ref{pro2-7}, a
contradiction. If $s=1$, then $|V(T_1)|=n-2$. We claim that
$\Delta(T_1)=n-3$. Assume, to the contrary, that $\Delta(T_1)\leq
n-4$. we have $cfc(G)=h(G)=cfc(T_1)\leq n-4$ by Proposition
\ref{pro2-7}, Theorems \ref{th4-1} and \ref{th4-2}, a contradiction.
Thus, $\Delta(T_1)=n-3$, which implies that $G=K_{1,n-1}+e$ or $G=U_n$.
The proof is complete.\qed
\end{pf}

\vskip0.3cm

Let $n$ be a natural number with $n\geq 7$. We now define a
sequence of graph classes, which will be used later.
\begin{itemize}
\item Let $U_n^1$ be a graph obtained from $U_{n-1}$ by adding a pendent
edge to a vertex of degree $2$ of $U_{n-1}$.

\item Let $U_n^2$ be a graph obtained from $U_{n-1}$ by adding a pendent edge
to a vertex of degree $1$ of $U_{n-1}$.

\item Let $U_n^3$ be a graph obtained from $K_{1,n-4}$ and $K_3$ by joining
a leaf vertex of $K_{1,n-4}$ and a vertex of $K_3$.

\item Let $U_n^4$ be a graph obtained from $K_{1,n-2}+e$ by adding a
pendent edge to a vertex of degree $2$ of $K_{1,n-2}+e$.

\item Let $U_n^5$ be a graph obtained from $K_{1,n-2}+e$ by adding a pendent
edge to a leaf vertex of $K_{1,n-2}+e$.

\item Let $U_n^6$ be obtained from $S_{4,n-4}$ by adding an edge joining the
two neighboring leaves of the vertex of degree $4$.

\item Let $W_n^1$ be a graph obtained from $C_4$ and $K_{1,n-4}$ by
identifying a vertex of $C_4$ and a leaf vertex of $K_{1,n-4}$.

\item Let $W_n^2$ be a graph obtained from $K_4$ and $K_{1,n-4}$ by
identifying a vertex of $K_4$ and a leaf vertex of $K_{1,n-4}$.

\item Let $W_n^3$ or $W_n^4$ be a graph obtained from $W_n^2$ by deleting an
edge of $K_4$.

\item Let $W_n^5$ be a graph obtained from $C_4$ and $K_{1,n-4}$ by
identifying a vertex of $C_4$ and the non-leaf vertex of
$K_{1,n-4}$.

\item Let $W_n^6$ be a graph obtained from $K_4$ and $K_{1,n-4}$ by
identifying a vertex of $K_4$ and the non-leaf vertex of
$K_{1,n-4}$.

\item Let $W_n^7$ or $W_n^8$ be a graph obtained from $W_n^6$ by deleting an
edge of $K_4$.
\end{itemize}

The following lemma is a preparation for the proof of Theorem
\ref{th4-8}.
\begin{lem}\label{lem4-7}
If $G\in \{U_n^1,\ldots, U_n^6,W_n^1,\ldots, W_n^8\}$ for $n\geq7$,
then $cfc(G)=n-4$.
\end{lem}
\begin{pf}
Note that each of these graph classes has the property that $h(G)\geq2$ and
there exists a unique component $T$ of $C(G)$ such that
$cfc(T)=h(G)$. Since $T=K_{1,n-4}$, it follows that
$cfc(G)=h(G)=n-4$ by Proposition \ref{pro2-7} and Theorem
\ref{th4-1}.
\end{pf}

\begin{thm}\label{th4-8}
Let $G$ be a connected graph of order $n\geq 5$. Then $cfc(G)=n-4$
if and only if $G$ is one of the following cases.
\begin{itemize}
\item[]$(i)$ $G$ is a tree with $\Delta(G)=n-4$ except for
$G=G_6$ in Fig. 1, where $ n\geq 7 $,

\item[]$(ii)$ $G$ is one of the 14 graph classes $\{U_n^1,\ldots, U_n^6,W_n^1,\ldots, W_n^8\}$ for $
n\geq 7 $,

\item[]$(iii)$ $G$ is a connected non-complete graph of order $6$
such that $G$ contains a cycle,

\item[]  $(iv)$ $G\in \{K_5, P_7, H_1, \cdots,
H_{5},H_{11}, H_{12}\}$, where $H_1, \ldots, H_{5}, H_{11}, H_{12}$
are shown in Fig. 3.
\end{itemize}
\end{thm}
\begin{pf}
If $G$ is a tree with $\Delta(G)=n-4$ except for $G=G_6$ in Fig. 1,
then $cfc(G)= n-4$ by Lemma \ref{lem2-2}, Theorems \ref{th4-1},
\ref{th4-2} and \ref{th4-6} for $n\geq 7$. If $G\in \{U_n^1,\ldots,
U_n^6,W_n^1,\ldots, W_n^8\}$ for $n\geq 7$, then $cfc(G)=n-4$ by
Lemma \ref{lem4-7}. Suppose that $G$ is a connected and non-complete
graph of order $6$. We first assume that $G$ is a tree. If
$\Delta(G)=2$, then $cfc(G)=3\neq n-4$, a contradiction. If
$\Delta(G)\geq 3$, then it follows from Lemma \ref{lem2-2} that
$cfc(G)\geq 3
> n-4$. Next, we deal with the case that $G$ has a cycle.
Let $C$ be the longest cycle. It is not hard to verify that
$cfc(G)=2$ in any case $|V(C)|=3, 4, 5,$ or $6$. It is clear that
$cfc(K_5)=cfc(P_7)=n-4$, and $cfc(H_i)=4=n-4$ for $1\leq i\leq 5$
and $cfc(H_{11})=5=n-4$ and $cfc(H_{12})=cfc(H_{13})=3=n-4$ in Fig.
3. Conversely, we let $G$ be a connected graph of order $n\geq 5$
with $cfc(G)=n-4$. If $n=5$, then $G=K_5$. If $n=6$, then we can
obtain that $G$ is neither a complete graph nor a tree. Next, we
assume $n\geq 7$. We distinguish the following two cases to show
this theorem.

\begin{figure}[h!]
\centering
\includegraphics[width=0.7\textwidth]{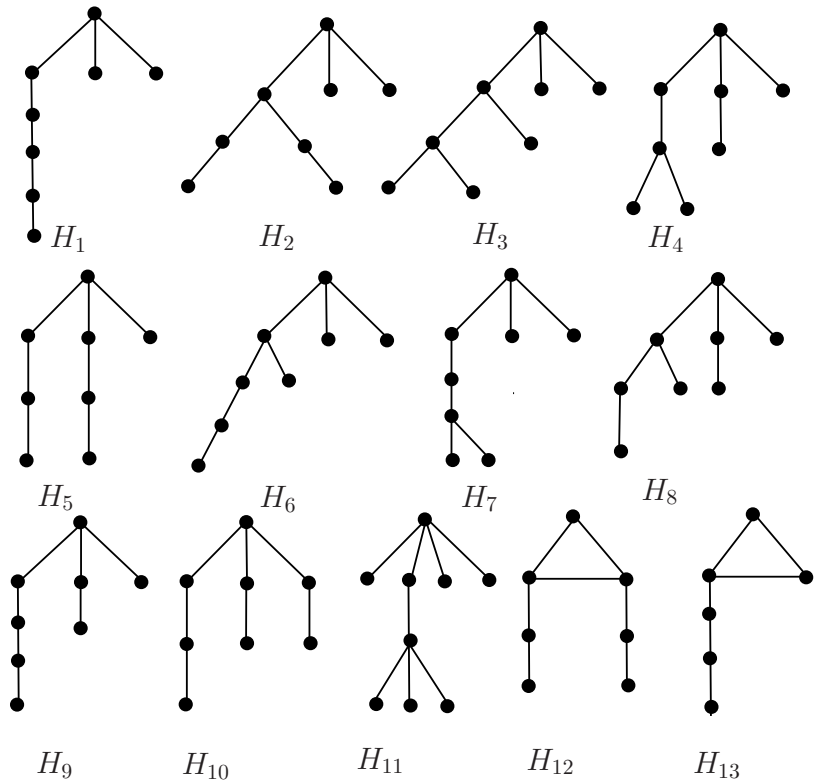}
\caption{Thirteen graphs in Theorem \ref{th4-8}}
\end{figure}

\begin{figure}[h!]
\centering
\includegraphics[width=0.7\textwidth]{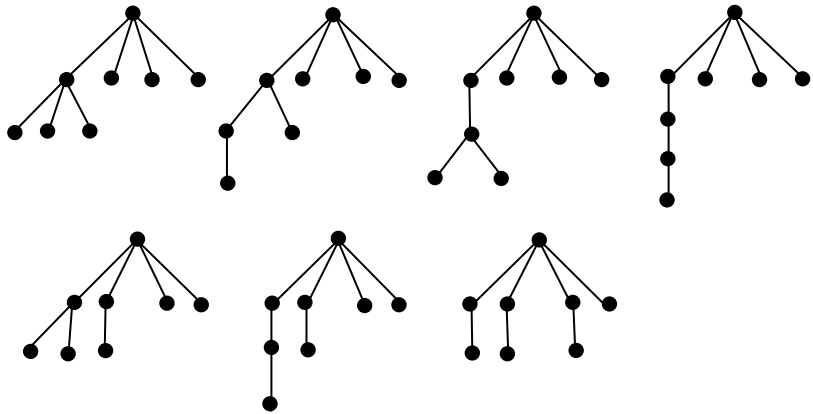}
\caption{Seven trees of order $8$ with maximum degree 4}
\end{figure}

\textbf{Case 1.} $G$ is a tree of order $n$ with $cfc(G)=n-4$.

If $n\geq 10$, then $\Delta(G)=n-4$ by Lemma \ref{lem3-2}. We only
need to consider the case $7\leq n\leq 9$. If $\Delta(G)=2$, then it
follows from Lemma \ref{lem2-1} that $G=P_7$. Suppose $\Delta(G)=3$.
If $n=7$, then it is easy to check that $G$ is not $G_6$ shown in Fig. 1.
If $n=8$, then $G$ is one of the first 10 trees shown in Fig. 3. One can
verify that $cfc(H_i)=4$ for $1\leq i\leq 5$ and $cfc(H_i)=3$ for $6\leq
i\leq 10$. Thus, $G \in \{H_1,\ldots,H_5\}$ shown in Fig. 3 in this case.
If $n=9$, then $cfc(G)\leq 4$, a contradiction. Suppose
$\Delta(G)=4$. If $n=7$, then it follows from Lemma \ref{lem2-2}
that $cfc(G)\geq \Delta(G)=4>n-4$, which is impossible. If $n=8$,
then $G$ is one of 7 trees shown in Fig. 4, and hence the result
clearly holds. If $n=9$, then $G=H_{11}$ shown in Fig. 3. Suppose
$\Delta(G)=5$. If $7\leq n\leq 8$, then $cfc(G)\geq \Delta(G)=5>n-4$
by Lemma \ref{lem2-2}, which is again impossible. If $n=9$, then it
is not hard to see that the result follows. If $\Delta(G)\geq 6$,
then it follows from Lemma \ref{lem2-2} that $cfc(G)\geq
\Delta(G)=6>n-4$, a contradiction.

\textbf{Case 2.} $G$ contains a cycle.

Since $G$ contains a cycle, it follows that there exists a
nontrivial block containing this cycle. Let $T_1,T_2,\ldots,T_s$ be
the components of $C(G)$ such that $h(G)=cfc(T_1)\geq cfc(T_2)\geq\cdots \geq cfc(T_s)$.
If $G$ has at least three nontrivial blocks,
then it follows from Lemma \ref{lem2-8} that
$cfc(G)\leq n-6$, which is impossible. Suppose $G$ has two
nontrivial blocks $B_1$ and $B_2$ with $|V(B_1)|\geq |V(B_2)|\geq
3$. It follows from the proof of Lemma \ref{lem2-8} that
$cfc(G)\leq h(G)+1\leq |V(T_1)| \leq
n-|V(B_1)|-|V(B_2)|+(2+s-1)-2(s-1)\leq n-|V(B_1)|-|V(B_2)|-s+3$.
If $|V(B_1)|\geq 4$, then $cfc(G)\leq
n+2-|V(B_1)|-|V(B_2)|\leq n-5$, a contradiction.
Suppose that $|V(B_1)|= |V(B_2)|=3$.
If $s\geq 2$, then $cfc(G)\leq n-|V(B_1)|-|V(B_2)|-s+3 \leq n-5$,
a contradiction. If $s=1$, then $cfc(G)=h(G)=cfc(T_1)\leq |V(T_1)|-1=n-4-1\leq n-5$
by Proposition \ref{pro2-7}, a contradiction. Next we only need to
consider that $G$ has only one nontrivial block $B$. If $|V(B)|\geq 6$,
then $cfc(G)\leq n+1-|V(B)|\leq n-5$ by Corollary \ref{cor2-9}, a
contradiction. Suppose $|V(B)|=5$. It follows from the proof of
Lemma \ref{lem2-8} that $cfc(G)\leq h(G)+1\leq |V(T_1)| \leq
n-|V(B)|+(1+s-1)-2(s-1)\leq n-|V(B)|-s+2\leq n-3-s$. If
$s\geq 2$, then $cfc(G)\leq n-5$, which is impossible. If $s=1$,
then $cfc(T_1)\leq |V(T_1)|-1 \leq n-4-1=n-5$, and hence
$cfc(G)=h(G)=cfc(T_1)\leq n-5$ by Proposition \ref{pro2-7}, which is
again impossible. Suppose $|V(B)|=4$. It follows from the proof of
Lemma \ref{lem2-8} that $cfc(G)\leq h(G)+1\leq |V(T_1)| \leq
n-|V(B)|+(1+s-1)-2(s-1)\leq n-|V(B)|-s+2\leq n-2-s$. If
$s\geq 3$, then $cfc(G)\leq n-5$, which is impossible. Suppose
$s=2$. If $|V(T_2)|\geq 3$, then $h(G)=cfc(T_1)\leq
|V(T_1)|-1 \leq n-|V(B)|-|V(T_2)|+2-1=n-6$, and so $cfc(G)\leq
h(G)+1=5$, a contradiction. If $|V(T_2)|=2$, then $|V(T_1)|\geq3$
and so $cfc(G)=h(G)=cfc(T_1)\leq |V(T_1)|-1\leq n-5$ by Proposition
\ref{pro2-7}, a contradiction. If $s=1$, then $|V(T_1)|=n-3$.
Note that $cfc(G)=h(G)=cfc(T_1)= n-4$ by Proposition \ref{pro2-7},
and hence it follows from Theorem \ref{th4-1} that $\Delta(T_1)=n-4$.
As a result, $G$ is a graph class in $\{W_n^1,\ldots,W_n^8\}$. Suppose
$|V(B)|=3$. Note that $s\leq 3$. Suppose $s=3$. As $n\geq7$, we have
$h(G)\geq 2$. If there exist two components $T_1, T_2$
such that $cfc(T_1)=cfc(T_2)=h(G)\geq 2$, then $cfc(G)\leq h(G)+1
\leq |V(T_1)|= n- |V(B)|-|V(T_2)|-|V(T_3)|+(1+3-1)\leq n-5$, a
contradiction. If there exists only one component $T_1$ such that
$cfc(T_1)=h(G)\geq 2$, then $cfc(G)=h(G)\leq |V(T_1)|-1= n-
|V(B)|-|V(T_2)|-|V(T_3)|+(1+3-1)-1\leq n-5$ by Proposition
\ref{pro2-7}, a contradiction. Suppose $s=2$. As $n\geq 7$ , we have
$h(G)\geq 2$. Without loss of generality, assume $cfc(T_1)\geq
cfc(T_2)$. If $cfc(T_1)= cfc(T_2)\geq 3$, then $cfc(G)\leq h(G)+1
\leq |V(T_1)|= n- |V(B)|-|V(T_2)|+(1+2-1)\leq n-5$, a contradiction.
If $cfc(T_1)=cfc(T_2)=2$, then the unique graph $H_{12}$ satisfies
$cfc(H_{12})=n-4$. If $cfc(T_1)>cfc(T_2)\geq 2$, then
$cfc(G)=h(G)\leq |V(T_1)|-1= n- |V(B)|-|V(T_2)|+(1+2-1)-1\leq n-5$
by Proposition \ref{pro2-7}, a contradiction. If
$cfc(T_1)>cfc(T_2)=1$, then $T_1$ is a star by Theorem \ref{th4-1},
and hence $G$ is $U_n^1$ or $U_n^4$. Suppose $s=1$. Since
$|V(T_1)|=n-2$, and $cfc(T_1)=h(G)=cfc(G)=n-4$ by Proposition
\ref{pro2-7}, it follows from Theorem \ref{th4-2} that $T_1$ is a
tree with $\Delta(T)=n-4$ for $n\geq6$ or $T=P_5$. Thus, $G\in
\{U_n^2,U_n^3,U_n^5,U_n^6\}$ or $G=H_{13}$. We complete the proof.\qed
\end{pf}

\section{Nordhaus-Gaddum-type theorem}

Note that if $G$ is a connected graph of order $n$,
then $m(G)\geq n-1$.
If both $G$ and $\overline{G}$ are connected,
then $n\geq 4$, since

$$2(n-1)\leq m(G)+m(\overline{G})=m(K_n)\leq\frac{n(n-1)}{2}.$$

In the sequel, we always assume that all graphs have
at least $4$ vertices, and both both $G$ and
$\overline{G}$ are connected.

\vskip0.3cm

\begin{thm}\label{th5-1}
Let $G$ and $\overline{G}$ be connected graphs of order
$n\ (n\geq 4)$. Then $4\leq cfc(G)+cfc(\overline{G})\leq n$
and $4\leq cfc(G)\cdot cfc(\overline{G})\leq2(n-2)$.
Moreover,
\begin{itemize}
\item[] $(i)$ $cfc(G)+cfc(\overline{G})=n$ or
$cfc(G)\cdot cfc(\overline{G})=2(n-2)$ if and only if one of $G$
    and $\overline{G}$ is a tree with maximum degree $n-2$
    or a $P_5$,
\item[] $(ii)$ the lower bounds are sharp.
\end{itemize}
\end{thm}

We proceed the proof of Theorem \ref{th5-1} by Proposition \ref{pro5-1},
Theorems \ref{th5-2} and \ref{th5-5}.

\vskip0.3cm

First, we study Nordhaus-Gaddum-type problem for conflict-free
connection number of graphs $G$ for the case both $G$ and $\overline{G}$
are $2$-edge-connected. It is known that if $G$ is a $2$-edge-connected
graph of order $4$, then $\overline{G}$ is disconnected. Thus, we need to
assume that $n\geq 5$ for the case that $G$ is a $2$-edge-connected in
the following.

\begin{pro}\label{pro5-1}
Let $G$ be a connected graph of order $n \ (n\geq 5)$. If both $G$
and $\overline{G}$ are $2$-edge-connected graphs,
then $cfc(G)+cfc(\overline{G})=4$ and $cfc(G)\cdot cfc(\overline{G})=4$.
\end{pro}
\begin{pf}
Since both $G$ and $\overline{G}$ are connected, it follows
that $G$ is neither an empty graph nor a complete graph.
Thus, $G$ and $\overline{G}$ are $2$-edge-connected and
noncomplete graphs, which implies $cfc(G)=2$ and
$cfc(\overline{G})=2$ by Corollary \ref{cor2-3}. The results hold. \qed
\end{pf}

\vskip0.3cm

Secondly, we turn to investigate the Nordhaus-Gaddum-type problem
for graphs $G$ such that $G$ is $2$-edge-connected and $\lambda(\overline{G})=1$.

\begin{thm}\label{th5-2}
Let $G$ be a $2$-edge-connected graph of order $n \ (n\geq 5)$ and $\lambda(\overline{G})=1$. Then $4\leq cfc(G)+cfc(\overline{G})\leq n$
and $4\leq cfc(G)\cdot cfc(\overline{G})\leq2(n-2)$.
Moreover,
\begin{itemize}
\item[] $(i)$ $cfc(G)+cfc(\overline{G})=n$ or
$cfc(G)\cdot cfc(\overline{G})=2(n-2)$
if and only if $\overline{G}=P_5$,
\item[] $(ii)$  the lower bounds are sharp.
\end{itemize}
\end{thm}
\begin{pf}
We proceed our proof by the following three claims.

\textbf{Claim 1.} $4\leq cfc(G)+cfc(\overline{G})\leq n$ and
$4\leq cfc(G)\cdot cfc(\overline{G})\leq2(n-2)$.

Since both $G$ and $\overline{G}$ are connected, it follows
that $G$ is neither an empty graph nor a complete graph.
Thus, $G$ is a $2$-edge-connected and noncomplete graph,
which implies that $cfc(G)=2$ and $cfc(\overline{G})\geq2$
by Corollary \ref{cor2-3}, Propositions \ref{pro2-11} and
\ref{pro3-1}. As a result, the lower bounds clearly hold.
It remains to verify the upper bound.
Firstly, we claim that $cfc(\overline{G})\leq n-2$.
Assume, to the contrary, that $cfc(\overline{G})=n-1$.
It follows that $\overline{G}=K_{1,n-1}$ by Theorem \ref{th4-1},
which implies that $G$ has an isolated vertex, a contradiction.
Thus, $cfc(\overline{G})\leq n-2$, and hence
$cfc(G)+cfc(\overline{G})\leq n$ and
$cfc(G)\cdot cfc(\overline{G})\leq2(n-2)$.\qed

\textbf{Claim 2.} $cfc(G)+cfc(\overline{G})=n$ or
$cfc(G)\cdot cfc(\overline{G})=2(n-2)$
if and only if $\overline{G}=P_5$.

\noindent\textbf{Proof of Claim 2:}
Since $cfc(G)=2$, it follows that 
$cfc(G)+cfc(\overline{G})= n$
or $cfc(G)\cdot cfc(\overline{G})=2(n-2)$
if and only if $cfc(\overline{G})=n-2$.
In order to complete our proof,
it is sufficient to show that
$cfc(\overline{G})=n-2$ if and only if $\overline{G}=P_5$.
Let $\overline{G}$ be a connected graph
such that $cfc(\overline{G})=n-2$.
It follows from Theorem \ref{th4-2}
that $\overline{G}$ is a tree with $\Delta(\overline{G})=n-2$
or $\overline{G}\in\{K_{1,3}+e, K_{2,2}, K_{2,2}+e, P_5\}$.
Suppose that $\overline{G}$ is a tree with $\Delta(\overline{G})=n-2$.
Let $v$ be a vertex of $\overline{G}$ with maximum degree $n-2$,
and $w$ be the unique vertex that is not adjacent to $v$ in $\overline{G}$.
Note that $vw$ is a cut-edge of $G$, which contradicts to that $G$ is $2$-edge-connected. If $\overline{G}\in\{K_{1,3}+e, K_{2,2}, K_{2,2}+e\}$,
then it is easy to see that $G$ is disconnected, a contradiction. Next, we only
need to consider the case $\overline{G}=P_5$. It is obtained that $G$ is a $2$-edge-connected and noncomplete graph of order $5$, and so $cfc(G)=2$.
\qed

\textbf{Claim 3.} The lower bounds in Theorem \ref{th5-2} are sharp.

\noindent\textbf{Proof of Claim 3:} The following example shows that the
lower bounds in Theorem \ref{th5-2} are best possible. Suppose that $H$ is
a $2$-edge-connected and noncomplete graph, and put one pair of nonadjacent
vertices $x,y$ of $H$. Let $u,v$ be two new vertices with $u,v \notin V(H)$, and 
let $\overline{G}$ be a graph such that $V(\overline{G})=V(H)\cup\{u,v\}$ and
$E(\overline{G})=E(H)\cup ux\cup uy\cup uv$. Clearly,
$\lambda(\overline{G})=1$ and $G$ is $2$-edge-connected.
It follows from Corollary \ref{cor2-3} and Lemma \ref{lem2-0} that
$cfc(G)=2$ and $cfc(\overline{G})=2$.\qed
\end{pf}

\vskip0.3cm

Finally, we discuss the Nordhaus-Gaddum-type problem for conflict-free connection number of graphs $G$ such that  $\lambda(G)=\lambda(\overline{G})=1$.

\vskip0.3cm

The following two lemmas are preparations for the proof of Theorem
\ref{th5-5}.
\begin{lem}\label{lem5-3}\cite{AH2}
A graph $G$ with $p$ points satisfies the condition $\lambda(G)=\lambda(\overline{G})=1$ if and only if
$G$ is a connected graph with a bridge and $\Delta=p-2$.
\end{lem}

\begin{lem}\label{lem5-4}
Let $G$ be a connected graph of order $n \ (n\geq 4)$.
If $\lambda(G)=1$ and $\lambda(\overline{G})=1$,
then at least one of $G$ and $\overline{G}$ has
conflict-free connection number $2$.
\end{lem}

\begin{pf}
By contradiction. Suppose that $cfc(G)\geq 3$ and
$cfc(\overline{G})\geq3$. Since $cfc(G)\geq 3$,
it follows from Lemma \ref{lem2-6} that $ h(G)\geq 2$.
At first, we present the follow claim.

\textbf{Claim 1.} $G$ has at least three cut-edges.

\noindent\textbf{Proof of Claim 1:} Let $T$ be the component
of $C(G)$ such that $ cfc(T)=h(G)\geq 2$. It follows that $m(T)\geq 2$.
Suppose $m(T)=2$. It follows that $h(G)=2$ and $cfc(G)=3$. By Propositon
\ref{pro2-7}, there exists another component $T'$ of $C(G)$ such that
$cfc(T')=h(G)=2$. Clearly, $m(T')\geq2$. Note that each edge of both $T$
and $T'$ is a cut-edge. Then $G$ has at least four cut-edges.
Suppose $m(T)\geq 3$. Since every edge of $T$ is a
cut-edge, it follows that $G$ has at least three cut-edges.\qed

It is obtained that $G$ has a cut-edge and a vertex $v$ of
degree $n-2$ by Lemma \ref{lem5-3}. Let $X=\{x_1,\ldots,x_{n-2}\}$ be
the neighborhood of $v$, and $w$ be the unique vertex
that is not adjacent to $v$. Note that $w$ is adjacent to at least one
vertex in $X$. Without loss of generality, assume that $wx_1\in E(G)$.
Let $T$ be a spanning tree induced on the set of edges 
$\{vx_1,\ldots,vx_{n-2},wx_1\}$ of $G$. Next, we give another claim.

\textbf{Claim 2.} If $e$ is a cut-edge of $G$, then $e\in E(T)$.

\noindent\textbf{Proof of Claim 2:} Assume, to the contrary, that
$e$ is a cut-edge of $G$ and $e\notin E(T)$. Since $T$ is a
spanning tree of $G-e$, it follows that $G-e$ is connected,
which is impossible.\qed

By Claims $1$ and $2$, there exist at least two cut-edges
incident with $v$ in $G$, say $vx_i$ and $vx_j$ for
$1\leq i< j\leq n-2$. If $i\geq2$, then $d(x_i)=d(x_j)=1$ and
$n\geq 5$. Noticing that $d_{\overline{G}}(x_i)=d_{\overline{G}}(x_j)=n-2$,
one can easily obtain that $\overline{G}[V\setminus v]$
contains $K_{2,n-3}$ as a subgraph $H_1$.
Let $H'$ be a spanning subgraph of $\overline{G}$ obtained
from $H_1$ by joining $w$ to $v$.
Then $cfc(H')=2$ by Lemma \ref{lem2-0}, which implies that 
$cfc(\overline{G})\leq cfc(H')=2$ by Proposition \ref{pro2-5},
a contradiction. Suppose $i=1$. Note that $wx_1$ is also a cut-edge,
and $d_{\overline{G}}(w)=d_{\overline{G}}(x_j)=n-2$.
If $n=4$, then $G=\overline{G}=P_4$, thus $cfc(G)=cfc(\overline{G})=2$,
a contradiction. If $n=5$, then it is easy to see that $cfc(\overline{G})=2$, a contradiction. If $n\geq6$, then $\overline{G}[V\setminus \{v,x_1\}]$ contains $K_{2,n-4}$ as a subgraph
$H_2$. Let $H''$ be a spanning subgraph of $\overline{G}$ obtained from
$H_2$ by adding two edges $wv$ and $x_jx_1$. Then $ cfc(H'')=2$ by Lemma
\ref{lem2-0}, which means that $cfc(\overline{G})\leq cfc(H'')=2$ by Proposition
\ref{pro2-5}, a contradiction.\qed
\end{pf}

\begin{thm}\label{th5-5}
Let $G$ be a connected graph of order $n\geq4$. If $\lambda(G)=1$
and $\lambda(\overline{G})=1$, then $4\leq cfc(G)+cfc(\overline{G})\leq n$
and $4\leq cfc(G)\cdot cfc(\overline{G})\leq2(n-2)$. Moreover,
\begin{itemize}
\item[] $(i)$ $cfc(G)+cfc(\overline{G})=n$ or
$cfc(G)\cdot cfc(\overline{G})=2(n-2)$ if and only if
one of $G$ and $\overline{G}$ is a tree with maximum degree
$n-2$.
\item[] $(ii)$ $cfc(G)+cfc(\overline{G})=4$ or $cfc(G)\cdot cfc(\overline{G})=4$
if and only if each of $G$ and $\overline{G}$ has the property that there exist at most two cut-edges incident with the vertex of maximum degree.
\end{itemize}
\end{thm}

\begin{pf}
We proceed our proof by the following three claims.

\textbf{Claim 1.} $4\leq cfc(G)+cfc(\overline{G})\leq n$ and
$4\leq cfc(G)\cdot cfc(\overline{G})\leq2(n-2)$.

\noindent\textbf{Proof of Claim 1:} Since both $G$ and $\overline{G}$ are
connected, it follows that $G$ is neither an empty graph nor a complete graph.
Thus, $cfc(G)\geq 2$ and $cfc(\overline{G})\geq2$ by Propositions
\ref{pro2-11} and \ref{pro3-1}. As a result, the lower bounds trivially
hold. For the upper bound, we first claim that $cfc(G)\leq n-2$.
Suppose, to the contrary, that $cfc(G)=n-1$. Then $G=K_{1,n-1}$
by Theorem \ref{th4-1}. It follows that $\overline{G}$ has an isolated
vertex, a contradiction. With a similar argument, one can show that
$cfc(\overline{G})\leq n-2$. Since at least one of $G$ and $\overline{G}$
has conflict-free connection number $2$ by Lemma \ref{lem5-4}, then
$ cfc(G)+cfc(\overline{G})\leq n$. \qed

\textbf{Claim 2.} $cfc(G)+cfc(\overline{G})=n$ or $cfc(G)\cdot cfc(\overline{G})=2(n-2)$ if and only if one of $G$ and $\overline{G}$
is a tree with maximum degree $n-2$.

\noindent\textbf{Proof of Claim 2:} By Claim $1$ and Lemma \ref{lem5-4},
we have that $ cfc(G)+cfc(\overline{G})=n $
or $cfc(G)\cdot cfc(\overline{G})=2(n-2)$
if and only if $ cfc(G)=n-2 $ and $cfc(\overline{G})=2$, or $cfc(G)=2$ and $cfc(\overline{G})=n-2$. By symmetry, we only need to consider one case, say 
$cfc(G)=n-2 $ and $cfc(\overline{G})=2$. Suppose that $cfc(G)=n-2$. It follows
from Theorem \ref{th4-2} that $G$ is a tree with $\Delta(G)=n-2$, or
$G\in\{K_{1,3}+e, K_{2,2}, K_{2,2}+e, P_5\}$. If $G$ is a tree with
$\Delta(G)=n-2$, then $cfc(\overline{G})=2$ by Lemma \ref{lem5-4}. If
$G\in\{K_{1,3}+e, K_{2,2}, K_{2,2}+e\}$, then it is easy to see that $\overline{G}$
is disconnected. Next, we only need to consider the case $G=P_5$. Note that
$\overline{G}$ is a $2$-edge-connected and noncomplete graph of order $5$,
which contradicts to the fact that $\lambda(\overline{G})=1$. Thus, $G$ is a
tree with $\Delta(G)=n-2$. \qed

\textbf{Claim 3.} $cfc(G)+cfc(\overline{G})=4$ or
$cfc(G)\cdot cfc(\overline{G})=4$
if and only if each of $G$ and $\overline{G}$ has the property that there exist at most two cut-edges incident with the vertex of maximum degree.

\noindent\textbf{Proof of Claim 3:} At first, we show that $cfc(G)+cfc(\overline{G})=4$
if and only if each of $G$ and $\overline{G}$ has the property that there exist
at most two cut-edges incident with the vertex of maximum degree. For the
necessity. Let $G$ be a connected graph of order $n\geq4$ such that
$\lambda(G)=1$ and $\lambda(\overline{G})=1$ and $cfc(G)+cfc(\overline{G})=4$.
Suppose that one of $G$ and $\overline{G}$, say $G$, has the oppose property
that there exist at least three cut-edges incident with the vertex of
maximum degree. In order to make $G$ conflict-free connected, these
cut-edges need to be assigned three distinct colors. Thus,
$cfc(G)\geq 3$, a contradiction.

For the sufficiency. If each of $G$ and $\overline{G}$ has the property
that there exist at most two cut-edges incident with the vertex of
maximum degree. By symmetry, we only need to show that $cfc(G)=2$.
Note that $G$ has a cut-edge and a vertex $v$ of degree $n-2$ by
Lemma \ref{lem5-3}. Let $X=\{x_1,\ldots,x_{n-2}\}$ be the neighborhood of $v$, and $w$ be the unique vertex that is not adjacent
to $v$. Note that $w$ is adjacent to at least one vertex in $X$.
Without loss of generality, assume that $wx_1\in E(G)$. Suppose
that there exists only one cut-edge $vx_j$ incident with $v$. If $wx_1$ is not a cut-edge, then let $F_1=G[V \setminus \{x_j\}]$ and $F'=F_1+vx_j$. Noticing that $F_1$ is $2$-edge-connected, it follows that $cfc(F')=2$ by Lemma \ref{lem2-0}. Clearly, $F'$ is a spanning subgraph of $G$, and hence $cfc(G)\leq cfc(F')=2$ by Proposition \ref{pro2-5}. Thus $cfc(G)=2$. If $wx_1$ is a cut-edge, then let $H_1=G[V \setminus \{x_j,w\}]$ and $H'=H_1+vx_j+x_1w$. It is
obvious that $H_1$ is $2$-edge-connected. Then we have that
$cfc(H')=2$ by Lemma \ref{lem2-0}. Noticing that $H'$ is a spanning
subgraph of $G$, then $cfc(G)\leq cfc(H')=2$ by Proposition \ref{pro2-5}. Thus $cfc(G)=2$. Suppose that there exist two cut-edges $vx_i$ and $vx_j$ incident with $v$ ($i<j$). If $wx_1$ is not a cut-edge, then let $F_2=G[V \setminus \{x_i,x_j,\}]$ and $F''=F_2+vx_i+vx_j$. Note that $F_2$ is $2$-edge-connected, and $F''$ has only one component $T$ of $C(G)$. Thus, $cfc(F'')=h(F'')=cfc(T)=2$ by Proposition \ref{pro2-7}. Clearly, $F''$ is a spanning subgraph of $G$, and hence $cfc(G)\leq cfc(F'')=2$ by Proposition \ref{pro2-5}. Thus $cfc(G)=2$. If $wx_1$ is a cut-edge,
then let $H_2=G[V \setminus \{x_i,x_j,w\}]$ and $H''=H_2+vx_i+vx_j+x_1w$.
Clearly, $H_2$ is $2$-edge-connected. And it is known that $H''$ has two
components $T_1$ and $T_2$ of $C(G)$. Without loss of generality,
assume that $T_1$ is the path $x_ivx_j$ (or $vx_iw$), and $T_2$ is the
edge $x_1w$ (or $vx_j$). Thus, $cfc(H'')=h(H'')=cfc(T_1)=2$ by
Proposition \ref{pro2-7}. Clearly, $H''$ is a spanning subgraph of $G$,
and hence $cfc(G)\leq cfc(H'')=2$ by Proposition \ref{pro2-5}.
Thus $cfc(G)=2$.

With a similar argument, one can obtain that $cfc(G)\cdot cfc(\overline{G})=4$
if and only if each of $G$ and $\overline{G}$ has the property that there exist at most two cut-edges incident with the vertex of maximum degree. \qed
\end{pf}


\begin{thebibliography}{1}

\bibitem{AH}
M. Aouchiche, P. Hansen, A survey of Nordhaus-Gaddum type relations, {\it Discrete Appl. Math.} {\bf 161(4-5)} (2013), 466--546.

\bibitem{AH2}
J. Akiyama, F. Harary, A graph and its complement with
specified properties $\uppercase\expandafter{\romannumeral1}$:
Connectivity, {\it J. Math. \& Math. Sci.} {\bf 2(2)} (1979),
223--228.

\bibitem{ALLZ}
E. Andrews, E. Laforge, C. Lumduanhom, P. Zhang, On proper-path
colorings in graphs, {\it J. Combin. Math. Combin. Comput.}
{\bf 97} (2016), 189--207.

\bibitem{BM}
J.A. Bondy, U.S.R. Murty, {\it Graph Theory}, GTM $244$, Springer,
$2008$.

\bibitem{BFGMMMT}
V. Borozan, S. Fujita, A. Gerek, C. Magnant, Y. Manoussakis,
L. Montero, Zs. Tuza, Proper connection of graphs,
{\it Discrete Math.} {\bf 312} (2012), 2550--2560.

\bibitem{CLL}
L. Chen, X. Li, H. Lian, Nordhaus-Gaddum-type theorem for
rainbow connection number of graphs,
{\it Graphs \& Combin.} {\bf 29} (2013), 1235--1247.

\bibitem{CJMZ}
G. Chartrand, G.L. Johns, K.A. McKeon, P. Zhang, Rainbow
connection in graphs, {\it Math. Bohem.} {\bf 133} (2008), 85--98.

\bibitem{CJV}
J. Czap, S. Jendrol', J. Valiska,  Conflict-free connection of
graphs, Accepted by {\it Discuss. Math. Graph Theory}.

\bibitem{CKP}
P. Cheilaris, B. Keszegh. D. P$\acute{a}$lv$\ddot{o}$igyi,
Unique-maximum and conflict-free coloring for hypergraphs and tree
graphs, {\it SIAM J. Discrete Math.} {\bf 27} (2013), 1775--1787.

\bibitem{CT}
P. Cheilaris, G. T$\acute{o}$th, Graph unique-maximum and
conflict-free colorings, {\it J. Discrete Algorithms} {\bf 9}
(2011), 241--251.

\bibitem{DGM}
S. Dantas, S. Gravier, F. Maffray, Extremal graphs for the list-coloring version of a theorem of Nordhaus and Gaddum, {\it Discrete Appl. Math.} {\bf 141} (2004), 93--101.

\bibitem{ELR}
G. Even, Z. Lotker, D. Ron, S. Smorodinsky, Conflict-free
coloring of simple geometic regions with applications to frequency
assignment in cellular networks, {\it SIAM J. Comput.} {\bf 33}
(2003), 94--136.

\bibitem{GLQ}
R. Gu, X. Li, Z. Qin, Proper connection number of random graphs,
{\it Theoret. Comput. Sci.} {\bf 609(2)} (2016), 336--343.

\bibitem{HH}
F. Harary, T. W. Haynes, Nordhaus-Gaddum inequalities for domination
in graphs, {\it Discret. Math.} {\bf 155} (1996), 99--105.

\bibitem{HR}
F. Harary, R. W. Robinson, The diameter of a graph and its complement,
{\it Amer. Math. Monthly} {\bf 92} (1985), 211--212.

\bibitem{HV}
A. Hellwig, L. Volkmann, The connectivity of a graph and its complement, {\it Discrete Appl. Math.} {\bf 156} (2008), 3325--3328.

\bibitem{LLZ}
E. Laforge, C. Lumduanhom, P. Zhang,  Characterizations of graphs
having large proper connection numbers, {\it Discuss. Math. Graph
Theory} {\bf 36(2)} (2016), 439--453.

\bibitem{LM}
X. Li, C. Magnant, Properly colored notions of connectivity--a
dynamic survey, {\it Theory \& Appl. Graphs} {\bf 0(1)} (2015), Art.
2.

\bibitem{LY}
X. Li, Y. Mao, Nordhaus-Gaddum type results for the generalized edge-connectivity of graphs, {\it Discrete Appl. Math.} {\bf 185} (2015),  102--112.

\bibitem{LSS}
X. Li, Y. Shi, Y. Sun, Rainbow connections of graphs: A survey,
{\it Graphs \& Combin.}  {\bf 29} (2013), 1--38.

\bibitem{LS}
X. Li, Y. Sun,  Rainbow Connections of Graphs, Springer Briefs
in Math., Springer, New York, 2012.

\bibitem{LWY}
X. Li, M. Wei, J. Yue, Proper connection number and connected
dominating sets, {\it Theoret. Comput. Sci.} {\bf 607} (2015),
480--487.

\bibitem{LWYA}
D. Li, B. Wu, X. Yang, X. An, Nordhaus-Gaddum-type theorem for Wiener index of gaphs when decomposing into three parts, {\it Discrete Appl. Math.} {\bf 159(15)} (2011), 1594--1600.

\bibitem{NG}
E. A. Nordhaus, J. W. Gauddum: On complementary graphs,
{\it Amer. Math. Monthly} {\bf 63} (1956), 175--177.

\bibitem{PT}
J. Pach, G. Tardos, Conflict-free colourings of graphs and
hypergraphs, {\it Comb. Probab. Comput.} {\bf 18} (2009), 819--834.

\bibitem{SDK}
E. Shan, C. Dang, L. Kang, A note on Nordhaus-Gaddum inequalities for domination, {\it Discrete Appl. Math.} {\bf 136} (2004), 83--85.

\bibitem{ZW}
L. Zhang, B. Wu, The Nordhaus-Gaddum type inequalities of some chemical indices, {\it MATCH Commun. Math. Comput. Chem.} {\bf 54} (2005), 189--194.

\end{thebibliography}
\end{document}